\documentclass{amsart2000}
\usepackage{graphicx,pstricks, amssymb,  amscd2000}
\usepackage[all]{xy}
\swapnumbers \newtheorem*{theorem*}{Theorem}
\newtheorem{theorem}{Theorem}[section]
\newtheorem{lemma}[theorem]{Lemma}
\newtheorem{problem}[theorem]{Problem}
\newtheorem{proposition}[theorem]{Proposition}
\newtheorem{corollary}[theorem]{Corollary}

\theoremstyle{definition} 
\newtheorem{definition}[theorem]{Definition}
\newtheorem{example}[theorem]{Example}

\theoremstyle{remark} 
\newtheorem{remark}[theorem]{Remark}

\numberwithin{equation}{section}
\newcommand{\firef}[1]{Figure~{\rm\ref{#1}}}
\newcommand{\thref}[1]{Theorem~{\rm\ref{#1}}}

\newcommand{\leref}[1]{Lemma~{\rm\ref{#1}}}
\newcommand{\coref}[1]{Corollary~{\rm\ref{#1}}}

\newcommand{\deref}[1]{Definition~{\rm\ref{#1}}}
\newcommand{\exref}[1]{Example~{\rm\ref{#1}}}

\newcommand{\reref}[1]{Remark~{\rm\ref{#1}}}

\newcommand{\seref}[1]{Section~{\rm\ref{#1}}}

\newcommand{\fig}[1]
{\raisebox{-0.5\height}%
  {\includegraphics{#1}}}

\newcommand{\newfig}[2]
{\raisebox{-0.5\height}%
{\includegraphics{#1}{#2}}
}

\newcommand{\st}{\; | \;}                     
\newcommand{\ttt}{\otimes}                    
\newcommand{\tta}{\otimes_A}                  
\newcommand{\tbox}{\boxtimes} 

\newcommand{\<}{\langle} 
\renewcommand{\>}{\rangle}
\newcommand{\injto}{\hookrightarrow}          
\newcommand{\isoto}{\xrightarrow{\sim}}       
\newcommand{\xxto}{\xrightarrow}              
\newcommand{\upcup}{\ar@{}@<-1ex>[u]^{\bigcup}} 

\newcommand{\one}{\mathbf{1}}     
\renewcommand{\i}{{\mathrm{i}}}   
\newcommand{\Cset}{\mathbb{C}}    
\newcommand{\Ctimes}{\mathbb{C}^\times}       
\newcommand{\Z}{\mathbb{Z}}       
\newcommand{\V}{{\mathcal{V}}}    
\newcommand{\C}{\mathcal{C}}      
\newcommand{\A}{\mathcal{A}}      
\newcommand{\F}{\mathcal{F}}      
\newcommand{\D}{(\rA)^G}          
\newcommand{\Do}{(\ro)^G}         
\newcommand{\DD}{\mathcal{D}}     
\renewcommand{\O}{\mathcal{O}}      
\renewcommand{\vec}{\mathcal{V}ec}

\newcommand{\al}{\alpha}

\newcommand{\de}{\delta} 

\newcommand{\la}{\lambda} 
\newcommand{\ph}{\varphi}
\newcommand{\Ph}{\Phi} 

\newcommand{\si}{\sigma} 
\renewcommand{\th}{\theta} 
\newcommand{\om}{\omega}
 \newcommand{\U}{U_q(\mathfrak{sl}_2)} 
\newcommand{\slthat}{\widehat{\mathfrak{sl}}_2}    

\DeclareMathOperator{\Res}{Res} 
\DeclareMathOperator{\Rep}{Rep}
\DeclareMathOperator{\Ind}{Ind} 
\DeclareMathOperator{\id}{id}
 
\DeclareMathOperator{\Hom}{Hom}
\DeclareMathOperator{\sym}{Sym}

\DeclareMathOperator{\tr}{tr} 
\DeclareMathOperator{\im}{Im}
\DeclareMathOperator{\Obj}{Obj}

\newcommand{\rA}{{\Rep A}}
\newcommand{\ro}{{\Rep^0 A}}

\begin{document}
\title{Modular  categories and orbifold models II} 
\author{Alexander Kirillov, Jr.}
\address{Department of Mathematics, SUNY at Stony Brook, 
         Stony Brook,  NY 11794, USA} 
\email{kirillov@math.sunysb.edu}
\urladdr{http://www.math.sunysb.edu/\textasciitilde kirillov/}
\thanks{The author was supported in part by NSF grant DMS9970473.}
                           

\maketitle
\section*{Introduction}
This paper is a continuation of the paper \cite{orb}. Its main goal is
to study so-called orbifold models in conformal field theory. In
mathematical language, this question can be reformulated as follows:
given a vertex operator algebra $\V$ and a finite group of
automorphisms $G$, describe the category $\C$ of modules over the
fixed-point algebra $\V^G$. In what follows, we assume that both $\V$
and $\V^G$ are rational VOA's and that categories of module over $\V,
\V^G$ are modular tensor categories.

It is well known that $\C$ can not be determined from the category of
$\V$-modules and the group $G$ only.  However, in the holomorphic case
(when $\V$ has a unique simple module, $\V$ itself) it was suggested
in \cite{DVVV}, \cite{DPR} and proved in \cite{orb} (for $\om=1$) that
if we additionally assume that $\C$ is modular, then $\C$ is
equivalent to the category of representations of a twisted Drinfeld
double $D^\om(G)$. (This result had also been announced --- without
full proof --- in \cite{mu}). Thus, in the holomorphic case the extra
data we need to completely determine $\C$ is the cocycle $\om \in
H^3(G, \Ctimes)$.

In this paper, we present some partial results in non-holomorphic
case. In particular, we show that the category $\C$ is completely
determined by the category of ``twisted'' $\V$-modules and the action
of the group $G$ on this category.  Our work is based on the results
of \cite{KO}, \cite{orb}, where it was shown that under certain
assumptions on $\V$, $\V^G$, this problem can be reformulated in the
language of tensor categories. We give this reformulation below, and
in the remainder of the paper only use the language of braided tensor
categories. Vertex operator algebras will not appear in the paper at
all.

Almost all results of this paper had been announced in \cite{mu};
however, \cite{mu} does not contain proofs of many important results,
referring the reader to a manuscript in preparation \cite{mu2}. Thus,
we feel that the present paper may help the readers to fill this gap.
Also, our methods are somewhat different from those of \cite{mu}.

Finally, it should be noted that some of the results of the current
paper had been proved in the language of VOA's in the paper \cite{DY};
see \reref{r:DY} for details. It also seems that there is a close
relation betweent he current paper and recent paper by Yamagami
\cite{Ya}. We plan to study this relation in detail in future
publications.

\section{Formulation of the problem}\label{s:formulation}
In this section, we list the main conventions used in this paper and
formulate the main problem in the language of tensor categories.

Throughout the paper, we keep the same notation as in \cite{KO},
\cite{orb}. In particular, 

\begin{itemize}
\item $\C$ is a semisimple rigid balanced braided
tensor category over $\Cset$ (later we will assume that $\C$ is modular), 

\item $A$ is a rigid commutative associative
algebra in $\C$ and $\th_A=\id$ (where $\th$ is the universal twist in
the category). 

\item $G$ a finite group acting faithfully by
automorphisms $\pi_g, g\in G$ on $A$ such that $A^G=\one$
\end{itemize}

As in   \cite{KO}, we use two categories of
$A$-modules, $\rA, \ro$ and functors $F\colon
\C\to\rA, G\colon \rA\to \C$. 

In such a situation, we will say that $\C$ is a ``$G$-orbifold of
$\ro$'' (this terminology is motivated by physics, as mentioned in the
introduction).

\begin{problem} Reconstruct $\C$ from the category $\ro$ and
  the group $G$ and probably some extra data.
\end{problem}

The paper \cite{orb} gave a proof of the result suggested in
\cite{DVVV}, namely, the answer to this question in the case when
$\ro\simeq \vec$ is the category of finite-dimensional vector
spaces (so-called holomorphic case). Here we will try to consider the
general case.

Before giving the answer, let us introduce some notation.
\begin{definition}
$\C_1\subset\C$ is the full subcategory generated (as an abelian
category) by simple objects $V_i$ such that $\<V_i, X\>\ne 0$ for
some   $X\in \ro$. 
\end{definition}

Then we have the following diagram:
\begin{equation}\label{e:maindiag}
\xymatrix{
\C \ar@<.5ex>[r]^F                  &            \rA \ar@<.5ex>[l]^G  \\
\C_1\ar@{}@<-1ex>[u]^{\bigcup}\ar@<.5ex>[r]^F  
         &  \ro \ar@{}@<-1ex>[u]^{\bigcup}\ar@<.5ex>[l]^G 
}
\end{equation}

Note that at the moment it is not clear why for every $V\in \C_1$,
$F(V)\in \ro$. This will be proved later (see \coref{c:F^0}).  

So, our goal can be formulated as follows: we need to reconstruct this
diagram from its lower right corner. However, we need more data. The
first piece of data is that in this situation, we have an ``action''
of $G$ on the category $\ro$.

As in  \cite[Theorem 4.7]{orb}), for  $X\in \rA, h\in G$ we
denote by $X^h$ the object of $\rA$ which coincides with $X$ as an
object of $\C$ but has the action of $A$ twisted by $\pi_h$. It
follows from the definition that $(X^h)^g=X^{hg}$. Thus, this
construction defines an ``action of $G$ by functors on $\rA$''. To be
precise, we have the following lemma.
\begin{lemma}\label{l:pig}
For $g\in G$ let $\Pi_g$ be a functor $\rA \to \rA$ defined
by 
\begin{equation}
\Pi_g(X)=X^{g^{-1}},
\end{equation}
and for $f\in \Hom_A(X,Y)$, $\Pi_g(f)=f$
considered as a morphism $X^{g^{-1}}\to Y^{g^{-1}}$. Then we have
canonical functor isomorphisms $\Pi_g\Pi_h=\Pi_{gh}$ and $\Pi_g(X\tta
Y)=\Pi_g(X)\tta \Pi_g(Y)$.
\end{lemma}
In particular, this shows that $G$ acts on
the set of isomorphism classes of objects in $\rA$ by 
\begin{equation}\label{e:gs}
g[X]=[X^{g^{-1}}].
\end{equation}

Note that we use $X^{g^{-1}}$ rather than $X^g$ in order to get the left
action; the formula $[X]\mapsto [X^g]$ would give right action.

The following construction allows us to pass from $\rA$ on which $G$
acts by functors to new category on which $G$ acts by automorphisms of
objects: 

\begin{definition} \label{d:catD}
Let $\D$ be the category with 
\begin{equation*}
\Obj\, \D=\{(X,\ph_X)\}\\
\end{equation*}
where $X$ is an object of $\rA$ and $\ph_X$ is a collection of
$\C$-morphisms $\ph_X(g)\colon X\to X, g\in G$ such that 
the following diagram is commutative:
\begin{equation}\label{e:phix}
\begin{CD}
A\ttt X               @>{\mu}>>           X\\
@V{\pi_g\ttt \ph_X(g)}VV             @V{\ph_X(g)}VV\\
A\ttt X               @>{\mu}>>           X
\end{CD}
\end{equation}
and 
\begin{equation}\label{e:ph}
\ph_X (g_1)\ph_X(g_2)=\ph_X(g_1g_2), \qquad \ph(1)=\id. 
\end{equation}

Morphisms in $\D$ are defined to be morphisms in $\rA$ which commute
with $\ph(g)$. 

Finally, we denote by $\Do$ the full subcategory in $\D$ with objects
$(X,\ph_X)$ such that $X\in \ro$. 
\end{definition}

Condition \eqref{e:phix} is equivalent to saying that
$\ph_X(g)$ is a morphism of $A$-modules $X\to X^g$, or, more
generally, for every $h$, $\ph_X(g)$ is a morphism of $A$-modules
$X^h\to X^{gh}$. 

Note that \eqref{e:ph} requires that $\ph(g)$ define an action of $G$
by $\C$-automorphisms of $X$. One might think that this is too
restrictive and we need to allow projective action. However, as we
explain later (see \thref{t:Dabelian}), it is not necessary.

Then the main result of this paper can be formulated as follows. 

\begin{theorem}
One has natural equivalences of braided tensor categories 
\begin{align*}
&\C\simeq \D,\\
&\C_1\simeq \Do.
\end{align*}
\end{theorem}

Proof of this theorem and further results will be given in
\seref{s:untwisted}, \seref{s:twisted1}. 

Comparing this with  the original question of recovering the
diagram~\eqref{e:maindiag} from its lower right corner, i.e. $\ro$, we
see that so far we have succeeded in reconstructing along the
horizontal lines: we can reconstruct $\C_1$ from $\ro$, and $\C$ from
$\rA$ (and the action of $G$ by automorphisms on $\rA$). 

\begin{remark}It should be noted that the bottom line of the diagram
  \eqref{e:maindiag} is ``modularization'' in the sense of \cite{B}.
\end{remark}

Reconstruction of $\C$ from $\C_1$ presents serious difficulties, some
of which will be discussed in \seref{s:conclusion}, and at the moment
seems to be out of reach   in non-holomorphic case.

\section{Model example}\label{s:model} 

Let $M$ be a finite group, and $N\subset M$ a normal subgroup. Let
$G=M/N$, $\C=\Rep M$. As in \cite{KO}, let $A=\F(M/N)$ be the algebra
of functions on $M/N=G$; we denote by $\de_{x}, x\in G$ the standard
basis of delta-functions in $A$. Then $A$  is a $\C$-algebra, and
$\rA=\ro$. Moreover, we have an action of $G$ by automorphisms on $A$:
\begin{equation}
\pi_g \de_x=\de_{xg^{-1}}, \quad x,g\in G.
\end{equation}
Thus we are in the situation discussed in \seref{s:formulation}, but
with the diagram~\eqref{e:maindiag} degenerating to a single line:
\begin{equation}\label{e:diag2}
\xymatrix{
\C \ar@<.5ex>[r]^F  &  \rA. \ar@<.5ex>[l]^G  
}
\end{equation}

In this situation, $\rA=\Rep N$, and the functors $F$
and $G$ are the restriction and induction functors, respectively (see
\cite{KO}), so the previous diagram becomes
\begin{equation}\label{e:diag2a}
\xymatrix{
\Rep M \ar@<.5ex>[r]^{\Res}  &  \Rep N. \ar@<.5ex>[l]^{\Ind}  
}
\end{equation}

Note that in particular, by \leref{l:pig}  we have an
action of $G$ by functor automorphisms on $\Rep N$. This action looks
especially simple if $M=G\ltimes N$ is a semidirect product: in this
case, $G$ acts on $N$ by conjugation, and the action on $\Rep N$ is
nothing by twisting the action of $N$ by this conjugation.

\begin{theorem}
One has canonical equivalence of categories 
$$
\C\simeq (\Rep N)^G.
$$
\end{theorem}
\begin{proof}
It follows from the definition that objects of $(\Rep N)^G=(\rA)^G$
are complex vector spaces $V$ with the following extra structure: 

 -- action of $M$

 -- $G$-grading: $V=\bigoplus_{g\in G} V_g$ such that $mV_g\subset
    V_{gm}, m\in M$

 -- action of $G$ which commutes with the action of $M$ and satisfies
    $\ph(g)V_x\subset V_{xg^{-1}}$

Define now the functors $(\Rep N)^G\to \Rep M$ and $\Rep M\to(\Rep
N)^G $ as follows: 

for $V\in (\Rep N)^G$, 
\begin{align*}
V\mapsto   V^G=& \text{collections of vectors } \{v_x\}_{x\in G}\\
               &\text{ such  that }v_x\in V_x, \ph(g)v_x=v_{xg^{-1}}
\end{align*}
and for $W\in \Rep M$, 
$$
W\mapsto W\ttt \F(G)
$$
which is considered as an object of $(\Rep N)^G$ with the grading
given by $(W\ttt \F(G))_x=W\ttt \Cset\de_x, x\in G$, action of $M$
given by $m(w\ttt \de_x)=mw\ttt \de_{mx}$ and action of $G$ by
$g(w\ttt \de_x)=w\ttt \de_{xg^{-1}}$.

It is easy to see that these functors are well-defined and inverse to
each other. 

\end{proof}

In this example, however, $\C$ is symmetric, so $\rA=\ro$. Let us
consider an example of a non-symmetric tensor category. 

Let $M$ be a finite group, $N\subset M$ --- a normal subgroup, and
$D(M)=\Cset[M]\ltimes\F(M)$ the Drinfeld double of $M$. Then we have
following commutative diagram of Hopf algebras:
\begin{equation}\label{e:diag3}
\xymatrix{
D(M) \ar[d]    \ar@{<-^{)}}[r]         &  \Cset[N]\ltimes \F(M)\ar[d]\\
\Cset[M]\ltimes \F(N)  \ar@{<-^{)}}[r] &  D(N)
}
\end{equation}
where the vertical arrows are restriction maps $\F(M)\to\F(N)$ and
horizontal arrows are inclusions $\Cset[N]\subset\Cset[M]$. 

As before, let $A=\F(G)=\F(M/N)$, which we consider as a module over
$D(M)$ by taking the grading to be identically $1$. One easily sees
that this is a commutative associative algebra in the category
$\C=\Rep D(M)$. As before, we have action of $G$ by automorphisms on
$A$. 
\begin{lemma}
In the situation above, one has canonical equivalences of categories 
\begin{equation}
\begin{aligned}
\rA&=\Rep (\Cset[N]\ltimes \F(M)),\\
\ro&=\Rep D(N),\\
\C_1&=\Rep (\Cset[M]\ltimes \F(N)).
\end{aligned}
\end{equation}
\end{lemma}
\begin{proof}
The proof is quite parallel to the calculations in \cite[Section 3]{orb}. 
\end{proof}

Thus, in this case the diagram \eqref{e:maindiag} is obtained 
by taking representations of the terms in diagram \eqref{e:diag3}:

\begin{equation}\label{e:diag4}
\xymatrix{
\Rep D(M) \ar@<.5ex>[r]^F    &  \Rep \Cset[N]\ltimes \F(M)\ar@<.5ex>[l]^G\\
\Rep \Cset[M]\ltimes \F(N) \upcup \ar@<.5ex>[r]^F &  
                  \Rep D(N)\upcup\ar@<.5ex>[l]^G
}
\end{equation}

\section{Category $\D$}\label{s:D}

In this section, we study some properties of the category $\D$
introduced in \deref{d:catD}. We start by giving examples of objects
in this category.

\begin{example}\label{ex:tx}
\begin{enumerate}
\item $X=A, \ph(g)=\pi_g$. From now on, we will denote this object of
  $\D$ by just $A$. 

\item Let $X\in \rA$. Let 
\begin{equation}
\Ind X=\bigoplus_{g\in G}X^g.
\end{equation}
Then $\Ind X$ has a natural action of $G$ by permutation of summands,
which shows that $\Ind X$ is an object of $\D$.
\end{enumerate}
\end{example}

Note that if $X\in \D$ and $V$ is a finite-dimensional representation
of $G$, then $V\tbox X\in \rA$ has a natural action of $G$ and thus a
structure of an object in $\D$ (see \cite[Section 1]{orb} for
definition of $\tbox$).  In other words, $\D$ is a module category
over $\Rep G$ of finite-dimensional complex representations of $G$. In
particular, applying this to $X=A$, we see that
\begin{equation}\label{e:v}
V\mapsto V\tbox A
\end{equation}
identifies $\Rep G$ with a subcategory in $\Do$. 

\begin{proposition}\label{p:dtensor}
  $\D$ has a natural structure of a rigid tensor category. $\Do$ has a
  natural structure of a balanced rigid braided tensor category, with
  the braiding inherited from $\C$.
\end{proposition}
\begin{proof}
Define $(X, \ph_X)\otimes (Y, \ph_Y)=(X\tta Y, \ph_X\ttt \ph_Y)$ and
$(X, \ph_X)^*=(X^*, (\ph_X(g^{-1}))^*)$, where for $f\in \Hom_\rA
(X,Y)$ we denote  by $f^*\in \Hom_\rA(Y^*, X^*)$ the adjoint
morphism. It is straightforward to check that so defined tensor product
and dual object satisfy all required properties.
\end{proof}

Note that the definition above also shows that the forgetful functor
$\D\to \rA$ is a tensor functor. 

\begin{remark} \label{r:dbraiding}
We can not define a braiding on $\D$ by just using the braiding in
$\C$: it will not be a morphism of $A$-modules (for this reason, $\rA$ in
general is not a braided category). However, we will show later (see
\thref{t:dbraiding}) that $\D$ does have a braided structure. 
\end{remark}

It is possible to give an explicit description of $\D$ as an abelian
category. Namely, let $S, S^0$ be the sets of isomorphism classes of
simple objects in $\rA$ and $\ro$ respectively. Let us fix for every
$s\in S$ a representative $X_s$. Formula \eqref{e:gs} defines an
action of $G$ on $S$. Let $\O\subset S$ be a $G$-orbit. Let $\F[\O]$
be the algebra of complex-valued functions on $\O$ and $\F^*[\O]$ the
group of non-vanishing functions on $\O$; both $\F[\O],\F^*[\O]$ are
naturally $G$-modules.  Then, as described in \cite{DY}, we have a
natural projective action of $G\ltimes \F[\O]$ by $A$-morphisms on
$X_\O=\bigoplus_{s\in \O}X_s$. Namely, the action of $\F[\O]$ is given
by $\de_s|_{X_s}=\id$, $\de_s|_{X_{s'}} =0$ for $s\ne s'$. To define
action of $G$, note that for every $s\in \O, g\in G$ there exists a
unique up to a constant $A$-isomorphism
\begin{equation}\label{e:phs}
\ph_s(g)\colon X_s^{g^{-1}}\isoto X_{g(s)}.
\end{equation}
Let us fix a choice of $\ph_s(g)$ and define $\ph_X(g)=\bigoplus_{s\in
  S} \ph_s(g)$. This gives a projective action of $G$ on $X_\O$.
Equivalently, we can say that $\O$ defines a cohomology class
$[\al]\in H^2(G, \F^*[\O])$ and we have a natural action of the
twisted algebra $\A_{\al}(G,\O)$ by $A$-morphisms on $X_\O$. We refer
the reader to \cite{DY} for details. 

In particular, if we choose $s\in
\O$ and denote
$$
G_s=\{g\in G\st gs=s\}=\{g\in G\st X_s^g\simeq X_s\},
$$
 then we have a projective action of $G_s$ on
$X_s$, or a true action of the twisted group algebra
$\Cset^{\al_s}[G_s]$.

\begin{remark}\label{r:leftright}
Please note that we consider $\F[\O]$ as a left $G$-module, whereas
\cite{DY} consider it as a right $G$-module. Thus, formulas of
\cite{DY} should be suitably modified for our setup (which is
actually more standard one). This modification is straightforward
enough and we leave it to the reader. 
\end{remark}

\begin{theorem}\label{t:Dabelian}
  As an abelian category,
  $$
  \D\simeq \bigoplus_\O \Rep\A_{\al^{-1}}(G,\O),
  $$
  where $\O$ runs through a set of orbits in $S$.
  The equivalence is given by
  \begin{align*}
  V&\mapsto V\tbox_{F[\O]}X_\O =\bigoplus_{s\in \O} V_s\tbox X_s, 
                   \qquad V\in \Rep  \A_{\al^{-1}}(G,\O)\\
  X&\mapsto V=\bigoplus_{s\in S} \Hom_A (X_s,X), \qquad X\in \D,
\end{align*}
where $X_\O=\bigoplus_{s\in \O} X_s$, $V=\bigoplus_{s\in \O} V_s$ is
the decomposition given by the action of $\F[\O]$, and the action of
$G$ on $X$ and $V$ is related by $\ph_X(g)=\bigoplus
\ph_{V_s}(g)\ttt\ph_s(g)$. Here $\ph_s(g)$ is as defined in
\eqref{e:phs} and $\ph_{V_s}(g)\colon V_s\to V_{gs}$ is the action of
$G$ on $V$.
\end{theorem}
\begin{proof}
  We need to check that the functors defined in the theorem are
  well-defined and inverse to each other (up to a functorial
  isomorphism), which is straightforward from the definition. The only
  step worth mentioning is that both $\ph_{V_s}$ and $\ph_s$ are
  projective representations of $G$, with cocycles $\al^{-1}$ and
  $\al$ respectively; thus, the tensor product $\ph_X(g)=\bigoplus
  \ph_{V_s}(g)\ttt\ph_s(g)$ is a true representation of $G$.
\end{proof}

This theorem shows that the structure of $\D$ as an abelian category
is completely determined by the set $S$ with the action of $G$ and the
cohomology class $[\al_S]=\bigoplus_{\O}[\al_\O]\in H^2(G, \F^*(S))$. 

Using the equivalence of categories 
$$
\Rep \Cset^{\al_s^{-1}}[G_s]\simeq \Rep \A_{\al^{-1}}(G,\O)
$$
(see \cite[Theorem 3.5]{DY}), the statement of the theorem can be
rewritten in simpler but less invariant way: there exists an
equivalence of abelian categories
$$
\D\simeq \bigoplus_s \Rep
  \Cset^{\al_s^{-1}}[G_s],
$$  
given by 
\begin{equation}\label{e:d->gs}
X\mapsto\bigoplus_s  \Hom_A(X_s, X)
\end{equation}
In both formulas, $s$ runs through a set of representatives of the orbits.

Combining this theorem with semisimplicity of $\Cset^\al[G_s]$ (see
\cite{K}), we get the following corollary.
\begin{corollary}\label{c:D}
\begin{enumerate}
\item 
$\D, \Do$ are semisimple.
\item
Simple objects in $\D$ \textup{(}respectively, $\Do$\textup{)} are 
\begin{equation}
X_{\la,\O}=V_\la \tbox_{\F[\O]} X_\O=\bigoplus_{s\in \O} (V_\la)_s\tbox X_s,
\end{equation}
where $\O$ is an orbit of the action of $G$ on the set $S$ of simple
objects in $\rA$ \textup{(}respectively, the set $S^0$ of simple objects in
$\ro$\textup{)}, and $V_\la$ is an irreducible representation of
$\A_{\al^{-1}}(G,\O)$.
\end{enumerate}
\end{corollary}

\begin{example}\label{ex:holom}
  Assume that $\ro\simeq \vec$, i.e. the only simple object in $\ro$
  is $A$. Then $V\mapsto V\tbox A$ defines an equivalence of tensor
  categories $\Rep G\simeq \Do$. This the so-called holomorphic case;
  it was studied in \cite{orb}. 
\end{example}

\begin{lemma}\label{l:tx}
Let $s\in S$ and let $\O$ be the orbit of $s$. Let $\Ind X_s$ be as in
\exref{ex:tx}. Then, as an object of $\D$, 
\begin{equation*}
\Ind X_s\simeq \bigoplus_{\la}V^*_\la \tbox X_{\la,\O},
\end{equation*}
where $\la$ runs over the set of irreducible representations of
$\Cset^{\al_s^{-1}}[G_s]$. In this formula, the multiplicity spaces
$V^*_\la\in \Rep \Cset^{\al_s}[G_s]$ are considered as vector spaces,
with no extra structure.
\end{lemma}
\begin{proof}
  It is immediate from the definitions that $\Ind X_s\simeq
  \bigoplus_{s'\in \O}\Cset[G_{s'}]\tbox X_{s'}$. Now the statement of
  the theorem 
  follows from the formula \eqref{e:d->gs} for the equivalence of
  categories $\D\simeq \Rep \Cset^{\al^{-1}_s}[G_s]$ and  identity
  $\Cset[G_s]\simeq V_\la^*\otimes V_\la$, where $\la$ runs over the
  set of irreducible representations of $\Cset^{\al_s^{-1}}[G_s]$. The
  latter formula holds for any finite-dimensional semisimple
  associative algebra. Note that if $V_\la$ is a representation of
  $\Cset^{\al_s^{-1}}[G_s]$, then $V_\la^*$ is naturally a
  representation of $\Cset^{\al_s}[G_s]$.
\end{proof}

\section{Untwisted sector}\label{s:untwisted}
Our first goal is to describe the ``untwisted sector'' of $\C$,
i.e. $\C_1$. Part of it has been done in \cite{orb} where we showed
that the subcategory in $\C$ generated by summands of $G(A)$ is
equivalent to $\Rep G$. 

In holomorphic case, this is the complete
description of the untwisted sector. In non-holomorphic case, more
work is needed.

Define a functor $\Ph\colon \D\to \C$ by 
$$
\Ph(X)=X^G,
$$ 
where we use the action of $G$ by $\C$-automorphisms of $X$ defined by
$\ph_X$. The definition of ``invariants'' $X^G$ is straightforward;
interested reader can find the details in \cite{orb}.

For example, it is easy to see that one has canonical isomorphisms 
\begin{align}
&\Ph(A)=\one, \label{e:Ph(A)}\\
&\Ph(\Ind X)=G(X), \qquad X\in \rA. \label{e:Ph(X)}
\end{align}

Note that $X^G$ is canonically a $\C$-sub-object in $X$; moreover, it
is easy to write a $\C$-morphism $\sym_X\colon X\to X$ which is a
projector on $X^g$:  $\sym^2=\sym, \im \sym=X^G$. It is given by
\begin{equation}\label{e:sym}
\sym=\frac{1}{|G|}\sum_{g\in G}\ph_X(g).
\end{equation}

\begin{theorem}\label{t:ut-main}
\begin{enumerate}
\item $\Ph$ is an equivalence of tensor categories $\D\simeq \C$
\item Restriction of $\Ph$ to $\Do\subset \D$ is an equivalence of
  braided tensor categories $\Do\simeq \C_1$, where $\C_1$ is the full
  subcategory in $\C$ generated \textup{(}as an abelian
  category\textup{)} by simple objects in $\C$ which appear in
  decomposition of $G(X), X\in \ro$.
\end{enumerate}
\end{theorem}
\begin{proof}
  
  The proof repeats, with suitable changes, the proof of Theorem~2.11
  in \cite{orb}, with $\Rep G$ replaced by $\D$. We sketch below those
  steps which are not completely identical.
  
  First, note that $\Ph(A)=\one$ and $\<\Ph(X_{\la,\O}), \one\>=0$ if
  $X_{\la,\O}$ is a simple object in $\D$ which is not isomorphic to
  $A$. Indeed, the first identity is obvious; the second follows from
  the fact that if $X_{\la,\O}\ne A$, then $\<G(X),
  \one\>_\C=\<X,A\>_\rA=0$.
  
  Next, define the functorial morphism $J\colon \Ph(X)\ttt \Ph(Y)\to
  \Ph(X\tta Y)$ as the following composition:
 
\begin{equation}\label{e:J}
X^G\ttt Y^G\injto (X\ttt Y)^G\to (X\tta Y)^G
\end{equation}
(the second morphism is induced by the canonical projection $X\ttt
Y\to X\tta Y$, see \cite{KO}).

Now we can repeat the same steps as in \cite{orb} to show that
$\Ph$ is compatible with associativity, commutativity (for $X,Y\in
\Do$) and unit isomorphisms, and with duality.  Then we show that the
subcategory $\DD\subset \D$ generated (as an abelian category) by
simple objects $X_{\la,\O}\in\D$ such that $\Ph(X_{\la,\O})\ne 0$ is
closed under tensor product and duality. By results of \cite{orb},
$\DD$ contains $\Rep G\subset \Do$. Now we can use the following lemma (cf.
\cite[Theorem 3.6]{orb}):
\begin{lemma}
  Let $\DD$ be a full subcategory in $\D$ which is closed under
  duality, tensor product, and taking sub-objects and contains $\Rep
  G$. Then $\DD$ is generated by $X_{\la, \O}$ where $\O$ runs over
  some set of orbits in $S$ and $\la$ runs over the set of all
  irreducible representations of $\A_\al(G,\O)$.
\end{lemma}

Assume that $\O$ is an orbit from the complement to the set of orbits
mentioned in the lemma, i.e. $\Ph(X_{\la,\O})=0$ for all $\la$. Choose
$s\in \O$ and consider $\Ind X_s$ as in \exref{ex:tx}. On one hand, by
\eqref{e:Ph(X)}, $\Ph((X_s)\tilde{})=X_s$. On the other hand, it follows from
\leref{l:tx} that $\Ph((X_s)\tilde{})=0$. This contradiction shows that
$\DD=\D$, and thus, $\Ph(X_{\la,\O})\ne 0$ for any simple $X_{\la,
  \O}$.

Now the same arguments as in \cite{orb} show that $\Ph$ is a tensor
functor which is an isomorphism on morphisms. Thus, $\Ph$ is an
equivalence $\D\simeq \C'$ where $\C'$ is a subcategory in $\C$
generated by $\Ph(X_{\la, \O)}$, each of which is a simple object in
$\C$. In other words, $\C'$ is the essential image of $\Ph$.

To show that $\C'=\C$, let $L$ be a simple object in $\C$. Then
$L\subset G(X)$ for some $X\in \rA$. Consider $\Ind X=\bigoplus X^g$.
Then $L\subset G(X)=\Ph (\Ind X)$ (see \eqref{e:Ph(X)}), which shows that
$L\in \C'$.

Replacing $\rA$ by $\ro$ we get part (2) of the theorem.
\end{proof}

\begin{remark}\label{r:DY}
This theorem contains as a special case the main result of \cite{DY},
which in our notation reads as follows: $\C_1$ is equivalent to $\Do$
as an abelian category. Note, however,  \cite{DY} does not discuss the
tensor structure of $\C_1$.
\end{remark}

\begin{theorem}\label{t:fg}
  Under the equivalence $\D\simeq \C$ described in \thref{t:ut-main},
  the functors $F\colon \C\to \rA, G\colon \rA\to \C$ are given by
\begin{align*}
&F((X, \ph_X))=X,\\ 
&G(X)=\Ind X.
\end{align*}
\end{theorem}
\begin{proof}
  For $G$, it easily follows from \eqref{e:Ph(X)}. For $F$, note that
  by definition  $F(X)=A\ttt X^G$. Define morphisms of
  $A$-modules $A\ttt X^G\to X, A\ttt X^G\to X$ by
\begin{align*}
&A\ttt X^G\injto A\ttt X\xxto{\mu}X\\
&X\xxto{i_A\ttt \id} A\ttt A\ttt X\xxto{\id\ttt\mu}A\ttt X 
   \to A\ttt X^G 
\end{align*}
where $X^G\injto X, X\to X^G$ are the canonical projection and
embedding and $i_A\colon \one\to A\ttt A$ is the rigidity morphism
(see \cite[Definition 1.9]{KO}).

We leave it to the reader to check that these two morphisms are
inverse to each other and thus, define isomorphism of $A$-modules
$X\simeq A\ttt X^G$. 
\end{proof}

\begin{corollary}\label{c:trans}
  $A$ is a ``transparent'', or ``central'', object in $\C_1$: for any
  $X\in \C_1$, $\check R_{XA}\check R_{AX}=\id_{A\ttt X}$.
\end{corollary}
Indeed, it suffices to note that $\Rep G\subset \Do$ is ``central'' in
$\Do$: for every $X\in \Do,Y\in \Rep G,\check R_{XY}\check
R_{YX}=\id$, which is obvious from the definitions.  In fact, it can
be shown (see \cite{mu}) that $\Do$ is exactly the centralizer in $\D$
of $\Rep G$, or, equivalently, $\C_1$ is the
centralizer in $\C$ of $\Rep G$.
\begin{corollary}\label{c:F^0}
For every $V\in \C_1$, $F(V)\in \ro$. 
\end{corollary}


\begin{remark}
The results of this section are almost identical to the results in
\cite{mu} if we note that the subcategory $\Rep G\subset \D\simeq \C$
is the same subcategory which is denoted by $S$ in \cite{mu}, and our
$\rA$ is the same as $\C\rtimes S$. The only difference is that
\cite{mu} uses unitarity of $\C$ in a non-trivial way. 
\end{remark}

\section{Twisted sector}\label{s:twisted1}
 {}From now on, let us assume that $\C$ is modular. 
Let us describe all of $\C$, not just the untwisted sector. 
Recall the definition of $g$-twisted $A$-module (see \cite[Definition
4.1]{orb}). 

\begin{theorem}\label{t:gtwisted}
  Every simple $X\in \rA$ is $g$-twisted for some $g\in G$.
\end{theorem}
\begin{proof}
  The proof is parallel to the proof of Theorem~4.3 in \cite{orb},
  with the following changes. Let $H$ be the following formal linear
  combination of objects in $\ro$:
\begin{equation}\label{e:H}
H=\bigoplus \frac{\dim_A X_s}{D^2}X_s,
\end{equation}
where $s$ runs over the set $S^0$ of simple objects in $\ro$, and
$D=\sqrt{\sum_s (\dim_A X_s)^2}$ (the normalization is chosen so that
$\dim_A H=1$). Then we have the following lemma:
\begin{lemma}
  If $X_i$ is a simple object in $\rA$, then
  $$
  \frac{1}{\dim A}
   \newfig{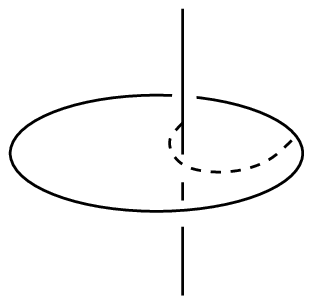}{
                      \rput(-1.4, 2.9){$X_i$}
                      \rput(-2.6, 1){$H$}
                      }              
                     =\de_{i0}\id_{X_i}
  $$
  where $i=0$ is the index of the unit object: $X_0=A$.
\end{lemma}
\begin{proof}[Proof of the lemma]
  First, by results of \cite[Lemma 5.4]{KO}, the left hand side is
  zero if $X_i\notin \ro$. If $X_i\in \ro$, then the result follows
  from the identity $(s^A)^2_{0i}=\de_{i0}$ (cf., for example,
  \cite[Corollary 3.1.11]{BK}).
\end{proof}

This lemma immediately implies the following result (cf.
\cite[Lemma~4.4]{orb})

\begin{corollary}
  $$
  \frac{1}{(\dim A)^2}\fig{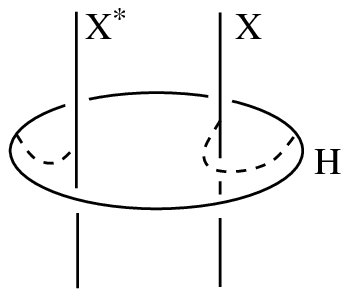}=\frac{1}{\dim
    X}\fig{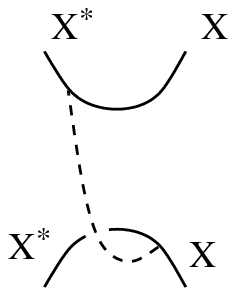}
  $$
\end{corollary}

Now we can repeat  the same steps as in the proof of Theorem~4.3 in
\cite{orb}, replacing as needed $A$ by $H$ and using \coref{c:trans}
instead of $\check R_{AA}^2=\id$.
\end{proof}

Let $\Rep_g A$ be the full subcategory in $\rA$ consisting of
$g$-twisted $A$-modules; then $\Rep_1 A$ is  the same category
which we had previously denoted $\ro$. It follows from
\thref{t:gtwisted} that
\begin{equation}\label{e:repg}
\rA=\bigoplus_{g\in G} \Rep_g A. 
\end{equation}
As in the holomorphic case (\cite[Theorem 4.7]{orb}), this grading
has some natural properties:
\begin{theorem}\label{t:repg}
 \begin{enumerate}
  \item If $X\in \Rep_g A$ then $X^h\in \Rep_{h^{-1}gh} A$.
  \item If $X_1\in \Rep_{g_1} A, X_2\in \Rep_{g_2} A$ then $X_1\tta
  X_2\in \Rep_{g_1g_2}A$.
  \item If $X\in \Rep_g A$, then $X^*\in \Rep_{g^{-1}}A$.
 \end{enumerate}
\end{theorem}
The proof of this theorem is completely parallel to the proof of
corresponding parts of \cite[Theorem 4.7]{orb} and is omitted.

Note that unlike holomorphic case, it is not true in general that
there is only one simple module in each of $\Rep_g A$, and it is not
true that every simple module  is invertible. 

We can now use decomposition \eqref{e:repg} to construct a braiding on
$\D$. For every $X\in \rA$, let  $\de_g\colon X\to X$ be the
projection on the $g$-twisted sector. i.e. $\de_g=\id$ for $X\in
\Rep_g A$ and $\de_g=0$ for $X\in \Rep_h A, h\ne g$.  
\begin{theorem}\label{t:dbraiding}
For $X, Y\in \D$, define  functorial morphism $\si_{X,Y}\in
\Hom_\C(X\ttt Y, Y\ttt X)$ by the following composition 
\begin{equation}\label{e:dbraiding}
X\ttt Y\xxto{\sum \de_g\ttt \ph_Y(g)} X\ttt Y\xxto{}Y\ttt X
\end{equation}
\textup{(}the second morphism is the usual commutativity isomorphism
in $\C$\textup{)}. Then $\si$ descends to a $\D$-morphism $X\tta Y\to
Y\tta X$ and defines a structure of a braided tensor category on $\D$.
\end{theorem}
\begin{proof}
  Let us start by showing that $\si$ descends to a $\C$-morphism
  $X\tta Y\to Y\tta X$. To do so, recall (see \cite{KO}) that $X\tta
  Y$ is defined as $X\ttt Y/I$, where $I=\im (\mu_1-\mu_2)$. Thus we
  need to show that $\si I\subset I$. 

  Without loss of generality, we can assume that $X\in \Rep_g A$. In
  this case, $\si=\id\ttt \ph_Y(g)$ and we can rewrite the composition
  $\si\circ (\mu_1-\mu_2)$ as shown in \firef{f:commut}, where the notation
  $f_1\equiv f_2$ stands for $\im(f_1-f_2)\in I$. Thus,
  $\si(\mu_1-\mu_2)\equiv 0$, or, equivalently, $\si(I)\subset I$.  
  
  The remaining parts of the theorem (i.e., that $\si$ is a morphism
  of $A$-modules, commutes with the action of $G$ and that it
  satisfies the hexagon axioms) are easily shown by direct calculation
  which we omit.
\begin{figure}[ht]
\begin{align*} 
&   \newfig{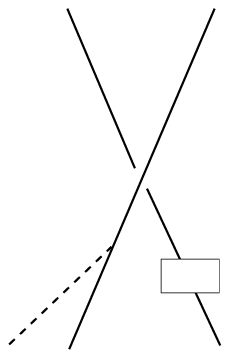}{
                   \rput(-0.8,1.15){$\scriptstyle\ph(g)$}
                   }
   -\newfig{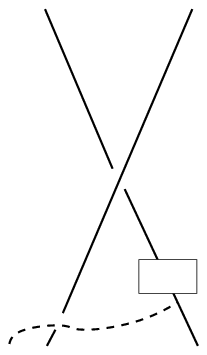}{
                   \rput(-0.8,1.1){$\scriptstyle\ph(g)$}
                     }
  = \newfig{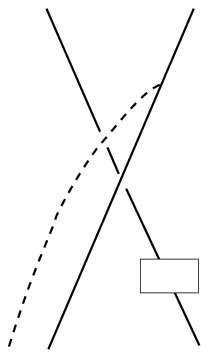}{
                   \rput(-0.8,1.15){$\scriptstyle\ph(g)$}
                   }
  -\newfig{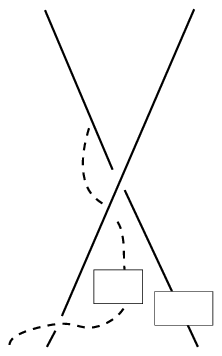}{
                     \rput(-0.75,0.78){$\scriptstyle\ph(g)$}
                     \rput(-1.4,0.95){${\scriptstyle\pi}_g$}
                    }\\
 &\qquad
   \equiv \newfig{commutc}{
                   \rput(-0.8,1.15){$\scriptstyle\ph(g)$}
                   }
  -\newfig{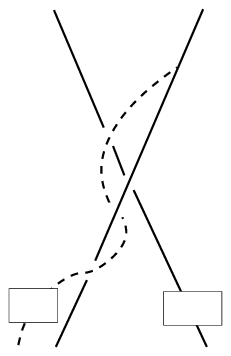}{                     
                     \rput(-0.75,0.78){$\scriptstyle\ph(g)$}
                     \rput(-2.35,0.8){${\scriptstyle\pi}_g$}
                    }
 \equiv 0
\end{align*}
\caption{Proof of \thref{t:dbraiding}}\label{f:commut}
\end{figure} 

\end{proof}

\begin{remark}
  It should be noted that $\sum \de_g\ttt g$ is the $R$-matrix for the
  Drinfeld double $D(G)$; thus, \eqref{e:dbraiding} combines the
  action of Drinfeld double on objects $X\in \D$ with the
  commutativity isomorphism in $\C$. For example, in the holomorphic
  case ($\ro=\vec$), as described in \cite{orb}, $\D\simeq\C$ is just
  the category of representations of (twisted) Drinfeld double
  $D^\om(G)$, and the commutativity isomorphism \eqref{e:dbraiding} is
  the usual commutativity isomorphism in $\Rep D^\om(G)$.
\end{remark}

\begin{remark}\label{r:gcrossed}
The decomposition \eqref{e:repg} and \thref{t:dbraiding} together are
equivalent to the statement that $\rA$ is a ``$G$-crossed braided
category'' as defined in \cite{Tu}, \cite{mu}. 
\end{remark}


\begin{proposition}
The functor $\Ph\colon \D\to \C$ identifies commutativity isomorphism
$\si$ in $\D$ defined by \eqref{e:dbraiding} with the commutativity
isomorphism in $\C$. 
\end{proposition}
\begin{proof}
Trivial: on $Y^G$, $\ph_Y(g)=\id$, so on $X^G\tta Y^G$,\\
 $\sum \de_g\ttt \ph_Y(g)=(\sum \de_g)\ttt \id=\id$.  
\end{proof}

\begin{corollary}
$\D$ is equivalent to $\C$ as a braided tensor category.
\end{corollary}

This result shows that the category $\C$ can be recovered from the
category $\rA$ and the action of $G$ on $\rA$. 

It would be desirable to  describe $\rA$ in terms of $\ro$
with some extra structure. However, at the moment we do not have such
a description. We do have some partial result, though. 

\begin{theorem}
If $\C$ is modular, then for all $g\in G$, $\Rep_g A\ne 0$.
\end{theorem}
\begin{proof}
Let $s$ be the $s$-matrix for $\C\simeq \D$. Let $X_i$ be a simple
object in $\D$ and $V_\la$ --- a simple object in $\Rep G\subset
\D$. Then it follows from the definition of $s$ and explicit formula 
\eqref{e:dbraiding} for the commutativity isomorphism in $\D$ that 
\begin{equation}\label{e:smatrix}
s_{V_\la, X_i}=\frac{1}{D_\C}\sum_g (\dim_A (X_i)_g)\tr_{V_\la}g 
\end{equation}
where $X=\bigoplus_g(X_i)_g$ is the decomposition of $X_i$ as an
$A$-module defined by \eqref{e:repg}. 

Now let $H=\{g\in G\st \Rep_g A\ne 0\}$. It immediately follows from
\thref{t:repg} that $H$ is a normal subgroup in $G$. Assume that $H\ne
G$. Then there exists a non-zero formal linear combination
$V=\bigoplus c_\la V_\la$ of irreducible representations of $G$ such
that $\tr_V g=0$ for all $g\in H$. Using \eqref{e:smatrix}, we get 
$$
s_{V, X_i}=\sum c_\la s_{V_\la, X_i}=0
$$
for all $X_i$, and thus the matrix $s$ is singular, which contradicts
the definition of a modular category. 
\end{proof}

\section{Conclusion}\label{s:conclusion}
Returning to the original question of recovering the
diagram~\eqref{e:maindiag} from its lower right corner, i.e. $\ro$, we
see that so far we have succeeded in reconstructing along the
horizontal lines: we can reconstruct $\C_1$ from $\ro$, and $\C$ from
$\rA$ (and the action of $G$ by automorphisms on $\rA$). 

However, we have failed to do the reconstruction along the vertical
lines: we can not recover $\C$ from $\C_1$, or $\rA$ from $\ro$. The
following example gives some insight in the difficulty of the problem:

\begin{example} 
  Let $\C$ be the semisimple subquotient of the category of
  representations of $\U$ at root of unity, $q=e^{\pi\i/l}$, as in
  \cite{KO}. Assume that $k=l-2$ is divisible by 4: $k=4m$. Then, as
  discussed in \cite{KO}, there is a unique structure of a commutative
  associative algebra on $A=\one+\de, \de=V_{k}$ (this corresponds to
  $D_{2m}$ in the ADE classification of ``subgroups'' in $\U$).
  Define an action of the group $\Z_2$ on $A$ by $\si|_{\one}=1,
  \si|_{\de}=-1$, where $\si$ is the non-trivial element of $\Z_2$.
  Then one easily sees that $\si$ is an automorphism of $A$ as a
  $\C$-algebra, and $A^{\si}=\one$, so this is an example of an
  $\Z_2$-orbifold. In this case, it follows from the results of
  \cite{KO} that $\C_1$ is the subcategory in $\C$ consisting of the
  modules with even highest weight (in physical terminology, integer
  spin). So in this case, recovering $\C$ from $\C_1$ is essentially
  equivalent to recovering the category of $\U$-modules from the
  subcategory of integer-spin modules. The classical analog of this
  problem would be recovering the category of representations of
  $\mathrm{SL}(2, \Cset)$ from the category of representations of
  $\mathrm{SL}(2, \Cset)/\Z_2=\mathrm{PSL}(2, \Cset)$. Even in the
  classical case,  we do not know any such construction. 
\end{example}

This example shows that a model question for recovering $\C$ from
$\C_1$ is reconstructing the category of representations of a group
$M$ knowing the normal subgroup $G\subset M$ (in the example above,
$G=\Z_2$) and the category $\Rep N, N=M/G$. This question is in a
sense dual to the model example discussed in \seref{s:model}, where
the roles of $G$ and $N$ were reversed.  It turns out that the dual
question is much harder. Of course, in principle this model question
is solvable: any finite group can be reconstructed from its category
of representations, so we recover $N$ from $\Rep N$, then use the
results of \seref{s:model}. Unfortunately, this recipe is very
indirect; even more importantly, it completely fails in the cases when
$\C_1$ is not a Tannakian category, so it is not a category of
representations of a group.



\end{document}